\documentclass[a4paper,11pt]{amsart}
\def\ladate{February 20, 2004. Final version, to appear. (v1: March 2002; 
v2: December 2002.)}

\def\From{From}
\newtheorem{theorem}{Theorem}[section]
\newtheorem{lemma}[theorem]{Lemma}
\newtheorem{corollary}[theorem]{Corollary}
\newtheorem{proposition}[theorem]{Proposition}

\newtheorem{note}{Note}

\theoremstyle{definition}
\newtheorem{definition}[note]{Definition}
\newtheorem{remark}[note]{Remark}

\newcommand{\NN}{{\mathbb N}}
\newcommand{\ZZ}{{\mathbb Z}}

\newcommand{\RR}{{\mathbb R}}
\newcommand{\CC}{{\mathbb C}}
\newcommand{\HH}{{\mathbb H}}

\newcommand{\cF}{{\mathcal F}}
\newcommand{\cG}{{\mathcal G}}
\newcommand{\cZ}{{\mathcal Z}}
\newcommand{\cA}{{\mathcal A}}
\newcommand{\cB}{{\mathcal B}}

\renewcommand{\Re}{{\rm Re}}
\renewcommand{\Im}{{\rm Im}}
\newcommand{\Res}{{\rm Res}}
\newcommand{\Un}{{\mathbf 1}}
\let\wh=\widehat
\let\wt=\widetilde

\begin{document}
\rightline{math.NT/0203120}
\rightline{\ }

\title[Two complete and minimal systems]{Two complete and
minimal systems associated with the zeros of the Riemann
zeta function}

\author{Jean-Fran\c cois Burnol}
\address{Universit\'e Lille 1, UFR de Math\'ematiques,
  Cit\'e scientifique M2, F-59655 Villeneuve d'Ascq, France}
\email{burnol@math.univ-lille1.fr}

\date{\ladate}



\begin{abstract}
We link together three themes which had remained separated
so far: the Hilbert space properties of the Riemann zeros,
the ``dual Poisson formula'' of Duffin-Weinberger
(also named by us co-Poisson formula), and the ``Sonine spaces'' of
entire functions defined and studied by de~Branges. We determine in which
(extended) Sonine spaces the zeros define a complete, or
minimal, system. We obtain some general results dealing with
the distribution of the zeros of the de-Branges-Sonine
entire functions. We draw attention onto some distributions
associated with the Fourier transform and which we
introduced in our earlier works.
\end{abstract}

\maketitle
\medskip

\section{The Duffin-Weinberger ``dualized'' Poisson formula (aka co-Poisson)}

We start with a description of the ``dualized Poisson
formula'' of Duffin and Weinberger (\cite{dufwein,
dufwein2}). We were not aware at the time of \cite{hab} that
the formula called by us \emph{co-Poisson formula} had been
discovered (much) earlier. Here is a (hopefully not too
inexact) brief historical account: the story starts with
Duffin who gave in an innovative 1945 paper \cite{duf1} a
certain formula constructing pairs of functions which are
reciprocal under the sine transform. As pointed out by
Duffin in the conclusion of his paper a special instance of
the formula leads to the functional equation of the
$L$-function $1 - \frac1{3^s} + \frac1{5^s} - \dots$ (as we
explain below, this goes both ways in fact). The co-Poisson
formula which we discuss later will stand in a similar
relation with the zeta function $1 + \frac1{2^s} +
\frac1{3^s}+\dots$, the pole of zeta adding its own special
touch to the matter. Weinberger extended in his dissertation
\cite{wein} this work of Duffin and also he found analogous
formulae involving Hankel transforms. Boas \cite{boas} gave
a formal argument allowing to derive Duffin type formulae
from the Poisson formula. However formal arguments might be
misleading and this is what happened here: formula
\cite[3.(iii)]{boas} which is derived with the help of a
purely formal argument looks like it is the \emph{co-Poisson
formula}, but is not in fact correct. It is only much later
in 1991 that Duffin and Weinberger \cite{dufwein} (see also
\cite{dufwein2}) published and proved the formula which, in
hindsight, we see now is the one to be associated with the
Riemann zeta function. They also explained its ``dual''
relation to the so-much-well-known Poisson summation
formula. In \cite{hab} we followed later a different
(esoterically adelic) path to the same result. As explained
in \cite{hab}, there are manifold ways to derive the
co-Poisson formula (this is why we use ``co-Poisson'' rather
than the ``dualized Poisson'' of Duffin and Weinberger). In
this Introduction we shall explain one such approach:
a re-examination of the Fourier meaning of the
functional equation of the Riemann zeta function.

When applied to functions which are compactly supported away
from the origin, the co-Poisson formula creates pairs of
cosine-tranform reciprocal functions with the intriguing
additional property that each one of the pair is constant in
some interval symmetrical around the origin. Imposing two
linear conditions we make these constants vanish, and this
leads us to a topic which has been invented by de~Branges as
an illustration, or challenge, to his general theory of
Hilbert spaces of entire functions (\cite{bra}),
apparently with the aim
to study the Gamma function, and ultimately also the Riemann
zeta function. The entire functions in these specific
de~Branges spaces are the Mellin transforms, with a Gamma
factor, of the functions with the vanishing property for
some general Hankel transform (the cosine or sine transforms
being special cases). These general ``Sonine Spaces'' were
introduced in \cite{bra64}, and further studied and
axiomatized by J. and V. Rovnyak in \cite{rov}. Sonine himself
never
dealt with such spaces, but in a study (\cite{son}) of
Bessel functions he constructed a pair of functions
vanishing in some interval around the origin and reciprocal
under some Hankel transform. An account of the Sonine spaces
is given in a final section of \cite{bra}, additional
results are to be found in \cite{bra92} and \cite{bra94}. As
the co-Poisson formula has not been available in these
studies, the way we have related the Riemann zeta function
to the Sonine spaces in \cite{hab} has brought a novel
element to these developments, a more intimate, and
explicit, web of connections between the Riemann zeta
function and the de~Branges spaces, and their extensions
allowing poles.

Although this paper is mostly
self-contained, we refer the reader to ``On
Fourier and Zeta(s)'' (\cite{hab}) for the motivating
framework and additional background and also to our Notes
\cite{cras2, cras3, cras4} for our results obtained so far
and whose aim is ultimately to reach
a better
understanding of some aspects of the Fourier Transform.

\vskip\the\baselineskip

Riemann sums $\sum_{n\geq1} F(n)$, or $\sum_{n\geq1} \frac1T
F(\frac nT)$, have special connections with, on one hand the
Riemann zeta function $\zeta(s) = \sum_{n\geq1} \frac1{n^s}$
(itself obtained as such a summation with $F(x) =
x^{-s}$), and, on the other hand, with the Fourier Transform.

In particular the functional equation of the Riemann zeta
function is known to be equivalent to the \emph{Poisson
summation formula}:
\begin{equation}\label{eq:poisson}
  \sum_{n\in\ZZ} \wt{\phi}(n) = \sum_{m\in\ZZ}{\phi}(m)
\end{equation}
which, for simplicity, we apply to a function $\phi(x)$ in
the Schwartz class of smooth quickly decreasing
functions. 

\begin{note}
We shall make
use of the following convention for the Fourier Transform:
$$\cF(\phi)(y) =
\wt{\phi}(y) = \int_\RR \phi(x)\,e^{2\pi i xy}dx$$
\end{note}

With a scaling-parameter $u\neq 0$,
\eqref{eq:poisson} leads to:
\begin{equation}\label{eq:poissonT}
  \sum_{n\in\ZZ} \wt{\phi}(nu) =
  \sum_{m\in\ZZ}\frac{1}{|u|}{\phi}(\frac mu)
\end{equation}
which, for the Gaussian $\phi(x) = \exp(-\pi x^2)$, gives
the Jacobi identity for the theta function (a function of
$u^2$). Riemann obtains from the theta identity  one of his
proofs of the functional equation of the zeta function,
which we recall here in its symmetrical form:
\begin{equation}\label{eq:functeq}
  \pi^{-s/2}\Gamma(\frac s2)\zeta(s) =
  \pi^{-(1-s)/2}\Gamma(\frac {1-s}2)\zeta(1-s)
\end{equation}

But there is more to be said on the Riemann sums
$\sum_{m\in\ZZ}\frac{1}{|u|}{\phi}(\frac mu)$ from the point
of view of their connections with the Fourier Transform than
just the Poisson summation formula \eqref{eq:poissonT};
there holds the \emph{co-Poisson intertwining formula}
(``dualized Poisson formula'' of Duffin-Weinberger
\cite{dufwein}), which reads:
\begin{equation}\label{eq:copoisson}
 \cF\left(\sum_{m\in\ZZ, m\neq0}\frac{g(m/u)}{|u|} -
  \int_\RR g(y)\,dy\right)(t) =
  \sum_{n\in\ZZ, n\neq0}\frac{g(t/n)}{|n|}
  - \int_\RR \frac{g(1/x)}{|x|}dx
\end{equation}

We show in \cite{hab} that it is enough to suppose for its
validity that the  integrals $\int_\RR \frac{g(1/x)}{|x|}dx$
and $\int_\RR g(y)\,dy$ are absolutely convergent. The
co-Poisson formula then computes the Fourier Transform of a
locally integrable function which is also tempered as a
distribution, the Fourier transform having the meaning given
to it by Schwartz's theory of tempered distributions. In the
case when $g(x)$ is smooth, compactly supported away from
$x=0$, then the identity is an identity of Schwartz
functions. It is a funny thing that the easiest manner to
prove for such a $g(x)$ that the sides of
\eqref{eq:copoisson} belong to the Schwartz class is to use
the Poisson formula \eqref{eq:poissonT} itself. So the
Poisson formula helps us in understanding the co-Poisson
sums, and the co-Poisson formula tells us things on the
Poisson-sums.

A most interesting case arises when the function $g(x)$ is
an integrable function, compactly supported away from $x=0$,
which turns out to have the property that the co-Poisson
formula is an identity in $L^2(\RR)$. The author has no
definite opinion on whether it is, or is not, an obvious
problem to decide which $g(x)$ (compactly supported away from
$x=0$) will be such that (one, hence) the two sides of the
co-Poisson identity are square-integrable. The only thing
one can say so far is that $g(x)$ has to be itself
square-integrable.

\begin{note}
Both the Poisson summation formulae \eqref{eq:poisson},
\eqref{eq:poissonT}, and the co-Poisson intertwining formula
\eqref{eq:copoisson} tell us $0=0$ when applied to
odd functions (taking derivatives leads to 
further identities which apply non-trivially to
odd-functions.) In all the following we deal only with even
functions on the real line. The square integrable among them
will be assigned squared-norm $\int_0^\infty
|f(t)|^2\,dt$. We let $K = L^2(0,\infty;dt)$, and we let
$\cF_+$ be the cosine transform on $K$:
$$\cF_+(f)(u) = 2\int_0^\infty \cos(2\pi tu)f(t)\,dt$$
The elements of $K$ are also tacitly viewed as even
functions on $\RR$.
\end{note}

Let us return to how the functional equation
\eqref{eq:functeq} relates with \eqref{eq:poissonT} and
\eqref{eq:copoisson}. The left-hand-side of
\eqref{eq:functeq} is, for $\Re(s)>1$, $\int_0^\infty
\sum_{n\geq1} 2e^{-\pi n^2 t^2}t^{s-1}\,dx$. An expression
which is valid in the critical strip is:
\[
0<\Re(s)<1\Rightarrow
\pi^{-s/2}\Gamma(\frac
s2)\zeta(s) = \int_0^\infty \left(\sum_{n\geq1} 2e^{-\pi n^2
t^2} - \frac1t\right)t^{s-1}\,dt
\]
More generally we have
the {M\"untz Formula} \cite[II.11]{tit}:
\begin{multline}\label{eq:zetaleft}
0<\Re(s)<1\Rightarrow\\
\zeta(s)\int_0^\infty \phi(t)
t^{s-1}\,dt = \int_0^\infty \left(\sum_{n\geq1} \phi(nt) -
\frac{\int_0^\infty \phi(y)dy}t\right)t^{s-1}\,dt
\end{multline}

We call the expression inside the parentheses the \emph{modified
  Poisson sum} (so the summation is accompanied with the
substracted integral). Replacing $\phi(t)$  with $g(1/t)/|t|$, with
$g(t)$ smooth, compactly supported away from $t=0$, gives
a formula involving a \emph{co-Poisson sum}:
\begin{multline}\label{eq:zetaright}
  0<\Re(s)<1\Rightarrow\\
  \zeta(s)\int_0^\infty g(t)
  t^{-s}\,dt = \int_0^\infty \left(\sum_{n\geq1}
  \frac{g(t/n)}n - {\int_0^\infty \frac{g(1/y)}y
  dy}\right)t^{-s}\,dt
\end{multline}

Let us now write $\wh{f}(s) = \int_0^\infty f(t)t^{-s}\,dt$
for the \emph{right} Mellin Transform, as opposed to the
\emph{left} Mellin Transform $\int_0^\infty
f(t)t^{s-1}\,dt$. These transforms are unitary
identifications of $K = L^2(0,\infty;dt)$ with
$L^2(s=\frac12+i\tau;d\tau/2\pi)$. Let $I$ be the unitary
operator $I(f)(t) = f(1/t)/|t|$. The composite $\cF_+\cdot
I$ is scale invariant hence diagonalized by the Mellin
Transform, and this gives, on the critical line:
\begin{equation}\label{eq:fourier}
\wh{\cF_+(f)}(s) = \chi(s)\wh{f}(1-s)
\end{equation}
with a certain function $\chi(s)$ which we obtain easily
from the choice $f(t) = 2\exp(-\pi t^2)$ to be
$\pi^{s-1/2}\Gamma(\frac {1-s}2)/\Gamma(\frac s2)$, hence
also $\chi(s) = \zeta(s)/\zeta(1-s)$.

The co-Poisson formula \eqref{eq:copoisson} follows then
from the functional equation in the form
\begin{equation}\label{eq:functeq2}
\zeta(s) = \chi(s)\zeta(1-s)
\end{equation}
together with \eqref{eq:fourier} and
\eqref{eq:zetaright}. And the Poisson formula
\eqref{eq:poissonT} similarly follows   from
\eqref{eq:functeq2} together with \eqref{eq:zetaleft}. We
refer the reader to \cite{hab} for further discussion and
perspectives.

The general idea of the equivalence between the Poisson
summation formula \eqref{eq:poissonT}  and the functional
equation \eqref{eq:functeq}, with an involvement of the
\emph{left Mellin Transform} $\int_0^\infty
f(t)t^{s-1}\,dt$, has been familiar and popular for many
decades. Recognizing that the \emph{right Mellin Transform}
$\int_0^\infty f(t)t^{-s}\,dt$ allows for a distinct
Fourier-theoretic interpretation of the functional equation
emerged only recently with our analysis \cite{hab} 
of the co-Poisson formula.

\section{Sonine spaces of de Branges and co-Poisson subspaces}

Let us now discuss some specific aspects of the co-Poisson
formula \eqref{eq:copoisson} (for an even function):
\begin{equation}\label{eq:copoisson2}
\cF_+\left(\sum_{m\geq1}\frac{g(m/t)}{|t|} -
  \wh{g}(0)\right)  =
\sum_{n\geq1}\frac{g(t/n)}{n} - \wh{g}(1)
\end{equation}
We are using the right Mellin transform  $\wh{g}(s) =
\int_0^\infty g(t)t^{-s}\,dt$. Let us take the (even)
integrable function $g(t)$ to be with its support in $[a,A]$
(and, as will be omitted from now on, also $[-A,-a]$ of
course), with $0<a<A$. Let us assume that the co-Poisson sum
$F(t)$ given by the right hand side belongs to $K =
L^2(0,\infty;dt)$. It has the property of being
equal to the constant $-\wh{g}(1)$ in $(0,a)$ and with its
Fourier (cosine) transform again constant in
$(0,1/A)$. After rescaling, we may always arrange that
$aA=1$, which we will assume henceforth, so that $1/A = a$
(hence, here, $0<a<1$).

So we are led to associate to each $a>0$ the sub-Hilbert
space $L_a$ of $K$ consisting of functions which are
constant in $(0,a)$ and with their cosine transform again
constant in $(0,a)$. Elementary arguments (such as the ones
used in \cite[Prop. 6.6]{hab}), prove that the $L_a$'s for
$0<a<\infty$ compose a strictly decreasing chain of
non-trivial infinite dimensional subspaces of $K$ with $K =
\overline{\cup_{a>0} L_a}$, $\{0\} = \cap_{a>0} L_a$, $L_a =
\overline{\cup_{b>a} L_b}$ (one may also show that
$\cup_{b>a} L_b$, while dense in $L_a$, is a proper
subspace).  This filtration is a slight variant on the
filtration of $K$ which is given by the \emph{Sonine spaces}
$K_a$, $a>0$, defined and studied by de~Branges in
\cite{bra64}. The Sonine space $K_a$ consists of the
functions in $K$ which are  vanishing identically, as well
as their Fourier (cosine) transforms, in $(0,a)$. The
terminology ``Sonine spaces'', from \cite{rov} and
\cite{bra}, includes spaces related to the Fourier sine
transform, and also to the Hankel transforms, and is used to
refer to some isometric spaces of analytic functions; we
will also call $K_a$ and $L_a$ ``Sonine spaces''. In the
present paper we use only the Fourier cosine tranform.

\begin{theorem}[De Branges \cite{bra64}]
Let $0<a<\infty$. Let $f(t)$ belong to $K_a$. Then its
completed right Mellin transform $M(f)(s)=
\pi^{-s/2}\Gamma(\frac s2)\wh{f}(s)$ is an entire
function. The evaluations at complex numbers $w\in\CC$ are
continuous linear forms on $K_a$.
\end{theorem}

We gave an elementary proof of this statement in
\cite{cras2}. See also \cite[Th\'eor\`eme 1]{cras3} for a
useful extension. Some slight change of variable is
necessary to recover the original de~Branges formulation, as
he ascribes to the real axis the r\^ole played here by the
critical line.  The point of view in \cite{bra64} is to
start with a direct characterization of the entire functions
$M(f)(s)$. Indeed a fascinating discovery of de~Branges is
that the space of functions $M(f)(s)$, $f\in K_a$ satisfies
all axioms of his general theory of Hilbert spaces of entire
functions \cite{bra} (we use the critical line where
\cite{bra} always has the real axis). It appears to be
useful not to focus exclusively on entire functions, and to
allow poles, perhaps only finitely many.

\begin{proposition}[{\cite[6.10]{hab}}]\label{thm:mellinLa}
Let $f(t)$ belong to $L_a$. Then its completed right Mellin
transform $M(f)(s)= \pi^{-s/2}\Gamma(\frac s2)\wh{f}(s)$ is
a meromorphic function in the entire complex plane, with at
most poles at $0$ and at $1$. The evaluations $f\mapsto
M(f)^{(k)}(w)$ for $w\neq0$, $w\neq1$, or $f\mapsto
\Res_{s=0}(M(f))$, $f\mapsto \Res_{s=1}(M(f))$ are
continuous linear forms on $L_a$. One has the functional
equations $M(\cF_+(f))(s) = M(f)(1-s)$.
\end{proposition}

We will write $Y^a_{w,k}$ for the vector in $L_a$ with
$$\forall f\in L_a\quad \int_0^\infty f(t)Y^a_{w,k}(t)\,dt =
M(f)^{(k)}(w)$$ This is for $w\neq 0,1$. For $w=0$ we have
$Y^a_0$ which computes the residue at $0$, and similarly
$Y^a_1$ for the residue at $1$. We are using the bilinear
forms $[f,g] = \int_0^\infty f(t)g(t)\,dt$ and not the
Hermitian scalar product $(f,g) = \int_0^\infty
f(t)\overline{g(t)}\,dt$ in order to ensure that the
dependency of $Y^a_{w,k}$ with respect to $w$ is analytic
and not anti-analytic. There are also evaluators $Z^a_{w,k}$
in the subspace $K_a$, which are (for $w\neq 0,1$)
orthogonal projections from $L_a$ to $K_a$ of the evaluators
$Y^a_{w,k}$.

\begin{definition}
We let $Y_a\subset L_a$ be the closed subspace of $L_a$
which is spanned by the vectors $Y^a_{\rho,k}$, $0\leq k
<m_\rho$, associated to the non-trivial  zeros $\rho$ of the
Riemann zeta function with multiplicity $m_\rho$.
\end{definition}

\begin{definition}
We let the ``co-Poisson subspace'' $P_a\subset L_a$, for
$0<a<1$, be the subspace of square-integrable functions
$F(t)$ which are co-Poisson sums of a function $g\in
L^1(a,A;dt)$ ($A = 1/a$).
\end{definition}

The subspace of $K_a$ defined analogously to $Y_a$ is denoted $Z_a$
(rather $Z_\lambda$) in \cite{hab}. The subspace of $K_a$
analogous to the co-Poisson subspace $P_a$ of $L_a$ is
denoted $W_a^\prime$ (rather $W_\lambda^\prime$) in
\cite{hab}. One has $W_a^\prime = P_a\cap K_a$. It may be
shown that if the integrable function $g(t)$, compactly
supported away from $t=0$, has its co-Poisson sum in $L_a$,
then $g$ is supported in $[a,A]$ and is square-integrable
itself.

\section{Statements of Completeness and Minimality}

It is a non-trivial fact that $P_a$ (and $W_a^\prime$ also,
as is proven in \cite{hab}) is closed. This is part of the
following two theorems.

\begin{theorem}\label{thmA}
The vectors $Y^a_{\rho,k}$, $0\leq k <m_\rho$, associated
with the non-trivial zeros of the Riemann zeta function, are
a minimal system in $L_a$ if and only if $a\leq1$. They are
a complete system if and only if $a\geq1$. For $a<1$ the
perpendicular complement to $Y_a$ is the co-Poisson subspace
$P_a$. For $a>1$ we may omit arbitrarily (finitely) many of
the $Y^a_{\rho,k}$'s and still have a complete system in
$L_a$.
\end{theorem}

\begin{theorem}\label{thmAprime}
The vectors $Z^a_{\rho,k}$, $0\leq k <m_\rho$, are a
minimal, but not complete, system for $a<1$. They are not
minimal for $a=1$, but the system obtained from omitting $2$
arbitrarily chosen among them (with the convention that one
either omits $Z^1_{\rho,m_\rho-1}$ and $Z^1_{\rho,m_\rho-2}$
or $Z^1_{\rho,m_\rho-1}$ and
$Z^1_{\rho^\prime,m_{\rho^\prime}-1}$) is again a minimal
system, which is also complete in $K_1$. In the case $a>1$
the vectors $Z^a_{\rho,k}$ are complete in $K_a$, even after
omitting arbitrarily (finitely) many among them.
\end{theorem}

\begin{remark}
This is to be contrasted with the fact that the evaluators
$Z^{1/\sqrt q}_{\rho,k}$ associated with the non-trivial
zeros of a Dirichlet $L$-function $L(s,\chi)$ (for an even primitive
character of conductor $q$) are a complete and
\emph{minimal} system in $K_{1/\sqrt q}$. Completeness was proven in
\cite[6.30]{hab}, and minimality is established as we will
do here for the Riemann zeta function.
\end{remark}

We use the terminology that an indexed collection of vectors
$(u_\alpha)$ in a Hilbert space $K$  is said to be minimal
if no $u_\alpha$ is in the closure of the linear span of the
$u_\beta$'s, $\beta\neq\alpha$, and is said to be complete
if the linear span of the $u_\alpha$'s is dense in $K$.  To
each minimal and complete system is associated a uniquely
determined dual system $(v_\alpha)$ with $(v_\beta,
u_\alpha) = \delta_{\beta\alpha}$ (actually in our $L_a$'s,
we use rather the bilinear form $[f,g] = \int_0^\infty
f(t)g(t)dt$). Such a dual system is necessarily minimal, but
by no means necessarily complete in general (as an example,
one may take $u_n = 1 - z^n$, $n\geq1$, in the Hardy space
of the unit disc. Then $v_m = -z^m$, for $m\geq1$, and they
are not complete).

For simplicity sake, let us assume that the zeros are all
simple. Then, once we know that $\zeta(s)/(s-\rho)$, for
$\rho$ a non-trivial zero, belongs to the space $\wh{L_1}$
of (right) Mellin transforms of elements of $L_1$, we then
identify the system dual to the $Y^{1}_{\rho,0}$'s, as
consisting of (the inverse Mellin transforms of) the
functions
$\zeta(s)/((s-\rho)\zeta^\prime(\rho)\pi^{-\rho/2}\Gamma(\rho/2))$.
Without any simplifying assumption, we still have that the
dual system is obtained from suitable linear combinations
(it does not seem very useful to spell them out explicitely)
of the functions $\zeta(s)/(s-\rho)^l$, $1\leq l \leq
m_\rho$, $\rho$ a non-trivial zero.

The proofs of \ref{thmA} and \ref{thmAprime} are a further
application of the technique of \cite[Chap.6]{hab}, which
uses a Theorem of Krein on Nevanlinna functions
\cite{krein,gorba}. Another technique is needed to establish the completeness
in $\wh{L_1}$ of the functions $\zeta(s)/(s-\rho)^l$, $1\leq
l \leq m_\rho$:

\begin{theorem}\label{thmB}
The functions $\zeta(s)/(s-\rho)^l$, for $\rho$ a
non-trivial zero and $1\leq l\leq m_\rho$ belong to
$\wh{L_1}$. They are minimal and complete in $\wh{L_1}$. The
dual system consists of vectors given for each $\rho$ by
triangular linear combinations of the evaluators
$Y^1_{\rho,k}$, $0\leq k <m_\rho$.
\end{theorem}

There appears in the proof of \ref{thmB} some computations
of residues which are reminiscent of a theorem of Ramanujan
which is mentioned in Titchmarsh \cite[IX.8.]{tit}.

The last section of the paper deals with the zeros of an
arbitrary Sonine functions, and with the properties of the
associated evaluators. We obtain in particular a density
result on the distribution of its zeros, with the help of
the powerful tools from the classical theory of entire
functions \cite{levin}.

\section{Aspects of Sonine functions}

\begin{note} We let $\wh{L_a}$ be the vector space of right
  Mellin transforms of elements of $L_a$ (and similarly for
  $\wh{K_a}$). They are square-integrable functions on the
  critical line, which, as we know from \ref{thm:mellinLa}
  are also meromorphic in the entire complex plane. We are
  not using here the Gamma-completed Mellin transform, but
  the bare Mellin transform $\wh{f}(s)$, which according to
  \ref{thm:mellinLa} has trivial zeros at $-2n$, $n>0$, and
  possibly a pole at $s=1$, and possibly does not vanish at
  $s=0$.
\end{note}

\begin{definition} We let  $\HH^2$ be the Hardy space of the
right half-plane $\Re(s) >\frac12$. We simultaneously view
$\HH^2$ as a subspace of $L^2(\Re(s)=\frac12,|ds|/2\pi)$ and
as a space of analytic functions in the right half-plane. We
also use self-explanatory notations such as $A^s \HH^2$.
\end{definition}

The right
Mellin transform is an isometric identification of
$L^2(1,\infty;dt)$ with $\HH^2$: this is one
of the famous theorems of Paley-Wiener \cite{palwin}, after
a change of variable. Hence, for $0<a$ and $A = 1/a$, the
right Mellin transform is an isometric identification of
$L^2(a,\infty;dt)$ with
$A^s\HH^2$. Furthermore, the right Mellin
transform is an isometric identification of
$\CC\cdot\Un_{0<t<a}+ L^2(a,\infty;dt)$ with
$\frac{s}{s-1}A^s\HH^2$. This leads to the
following characterization of $\wh{L_a}$:

\begin{proposition}\label{prop:La}
The subspace $\wh{L_a}$ of $L^2(\Re(s)=\frac12,|ds|/2\pi)$
consists of the measurable functions $F(s)$ on the critical
line which belong to $\frac{s}{s-1}A^s\HH^2$
and are such that $\chi(s)F(1-s)$ also belongs to
$\frac{s}{s-1}A^s\HH^2$. Such a function
$F(s)$ is the restriction to the critical line of an
analytic function, meromorphic in the entire complex plane
with at most a pole at $s=1$, and with trivial zeros at $s =
-2n$, $n\in\NN, n>0$.
\end{proposition}

\begin{proof}
We know already from \ref{thm:mellinLa} that functions
in $\wh{L_a}$ have the stated properties. If a function
$F(s)$ belongs to $\frac{s}{s-1}A^s\HH^2$,
viewed as a space of (equivalence classes of) measurable
functions on the critical line, then it is square-integrable
and is the Mellin transform of an element $f(t)$ of
$\CC\cdot\Un_{0<t<a}+ L^2(a,\infty;dt)$. We know that
the Fourier cosine transform of $f$ has $\chi(s)F(1-s)$ as
Mellin transform, so the second condition on $F$ tells us
that $f$ belongs to $L_a$.
\end{proof}

We recall that $\chi(s)$ is the function (expressible in
terms of the Gamma function) which is involved in the
functional equation of the Riemann zeta function
\eqref{eq:functeq2}, and is in fact the spectral multiplier
of the scale invariant operator $\cF_+\cdot I$, for the
right Mellin transform. 

\begin{note}
Abusively, we will say that
$\chi(s)F(1-s)$ is the Fourier transform of $F(s)$, and will
sometimes even write $\cF_+(F)(s)$ instead of
$\chi(s)F(1-s)$. It is useful to take note that if we write
$F(s) = \zeta(s)\theta(s)$ we then have $\chi(s)F(1-s) =
\zeta(s)\theta(1-s)$.
\end{note}

\begin{proposition}\label{zetaprop} The functions
$\zeta(s)/(s-\rho)^l$, $1\leq l \leq m_\rho$ associated with
the non-trivial zeros of the Riemann zeta function belong to
$\wh{L_1}$. \end{proposition}

\begin{proof} The function $F(s) = \zeta(s)/(s-\rho)^l$ is
square-integrable on the critical line. And $\chi(s)F(1-s) =
(-1)^l \zeta(s)/(s-(1-\rho))^l$. So we only need to prove
that $\frac{s-1}sF(s) = \frac{s-1}s\zeta(s)/(s-\rho)^l$
belongs to $\HH^2$. This is well-known to be
true of $\frac{s-1}s\zeta(s)/s$ (from the formula
$\zeta(s)/s =1/(s-1) -\int_1^\infty \frac{\{t\}}t t^{-s}dt$,
valid for $0<\Re(s)$), hence it holds also for
$\frac{s-1}s\zeta(s)/s^l$.  If we exclude a neigborhood of
$\rho$ then $s^l/(s-\rho)^l$ is bounded, so going back to
the definition of $\HH^2$ as a space of
analytic functions in the right half-plane with a uniform
bound of their $L^2$ norms on vertical lines we obtain the
desired conclusion.
\end{proof}

The following will be useful later:

\begin{proposition}\label{zeroprop}
If $G(s)$ belongs to $\wh{L_a}$ and 
$s(s-1)\pi^{-s/2}\Gamma(\frac s2)G(s)$ vanishes at $s=w$
then $G(s)/(s-w)$ again belongs to $\wh{L_a}$. If $G(s)$
belongs to $\wh{K_a}$ and $\pi^{-s/2}\Gamma(\frac s2)G(s)$
vanishes at $s=w$ then $G(s)/(s-w)$ again belongs to
$\wh{K_a}$.
\end{proposition}

\begin{proof}
We could prove this in the ``$t$-picture'', but will do it
in the ``$s$-picture''. We see as in the preceding proof
that $G(s)/(s-w)$ still belongs to
$\frac{s}{s-1}A^s\HH^2$. The entire function
$s(s-1)\pi^{-s/2}\Gamma(\frac s2)G(s)$ vanishes at $s=w$ so
$s(s-1)\pi^{-s/2}\Gamma(\frac s2)\cF_+(G)(s)$ vanishes at
$s=1-w$ and the same argument then shows that
$\chi(s)G(1-s)/(1-s-w)$ belongs to
$\frac{s}{s-1}A^s\HH^2$. We then apply Proposition
\ref{prop:La}. The statement for $K_a$ is proven analogously.
\end{proof}

\begin{note}
It is a general truth in all de~Branges' spaces that such a
statement holds for zeros $w$ off the symmetry axis (which
is here the critical line). This is, in fact, almost one of
the axioms for de~Branges' spaces. The possibility to divide
by $(s-w)$ if $w$ is on the symmetry axis depends on whether
the structure function $E$ (on this, we refer to \cite{bra})
is not vanishing or vanishing at $w$. For the Sonine spaces,
the proposition \ref{zeroprop} proves that the structure
functions $E_a(z)$ have no zeros on the symmetry axis. For
more on the $E_a(z)$'s and allied functions, see
\cite{cras3} and \cite{cras4}.
\end{note}

A variant on this gives:

\begin{lemma} If $F(s)$ belongs to $\wh{K_a}$ then $F(s)/s$
belongs to $\wh{L_a}$. \end{lemma}

\begin{proof} The function $F(s)/s$ (which is regular at
$s=0$) belongs to the space $A^s\HH^2$, simply from $1/|s| =
O(1)$ on $\Re(s)\geq\frac12$. Its image under the Fourier
transform is $\cF_+(F)(s)/(1-s)$ which belongs to
$\frac{s}{s-1}A^s\HH^2$.
\end{proof}

\begin{proposition}\label{dimlemma} One has
$\dim(L_a/K_a) = 2$. \end{proposition}

\begin{proof}
This is equivalent to the fact that the residue-evaluators
$Y^a_0$ and $Y^a_1$ are linearly independent in $L_a$, which
may be established in a number of elementary ways; we give
two proofs. Evaluators off the symmetry axis are always
non-trivial in de~Branges spaces so there is $F(s)\in
\wh{K_a}$ with $F^\prime(0)\neq0$ (one knows further \From\
\cite[Th\'eor\`eme 2.3.]{cras2} that any finite system of
vectors $Z^a_{w,k}$ in $K_a$ is a linearly independent
system).  So we have $F(s)/s = G(s)\in \wh{L_a}$ not
vanishing at $0$ but with no pole at $1$. Its ``Fourier
transform'' $\chi(s)G(1-s)$ vanishes at $0$ but has a pole
at $1$. This proves $\dim(L_a/K_a)\geq2$ and the reverse
equality follows from the fact that the subspace $K_a$ is
defined by two linear conditions.

For the second proof we go back to the argument of
\cite{cras2} which identifies the perpendicular complement
to $K_a$ in $L^2(0,\infty;dt)$ to be the closed space
$L^2(0,a)+\cF_+(L^2(0,a))$.  It is clear that $L_a$ is the
perpendicular complement to the (two dimensions) smaller
space $(L^2(0,a)\cap \Un_{0<t<a}^\perp)+\cF_+(L^2(0,a)\cap
\Un_{0<t<a}^\perp)$ and this proves \ref{dimlemma}.
\end{proof}

The technique of the second proof has the additional benefit:

\begin{proposition}\label{denselemma} The union
$\bigcup_{b>a} K_b$ is dense in $K_a$, and
$\bigcup_{b>a} L_b$ is dense in
$L_a$. \end{proposition}

\begin{proof} Generally speaking $\bigcap_{b>a}
(A_b + B_b) = (\cap_{b>a} A_b) + (\cap_{b>a} B_b)$ when we
have vector spaces indexed by $b>a$  with $A_{b_1}\subset
A_{b_2}$ and $B_{b_1}\subset B_{b_2}$ for $b_1<b_2$ and
$A_b \cap B_b = \{0\}$ for $b>a$. We apply this to $A_b =
L^2(0,b;dt)$ and $B_b = \cF_+(L^2(0,b;dt))$, as $K_a =
(A_a+B_a)^\perp$, $A_a = \cap_{b>a} A_b$,  $B_a = \cap_{b>a}
B_b$, and $A_a+ B_a$ is closed as a subspace of
$L^2(0,\infty;dt)$.
\end{proof}

\begin{proposition} 
The vector space $\bigcup_{b>a} K_b$ is properly included in
$K_a$ and the same holds for the respective subspaces of
Fourier invariant, or skew, functions (and similarly for
$L_a$).
\end{proposition}

\begin{proof}
Let $g\in K_b$, with $b>a$ and $g$ having the leftmost point
of its support at $b$. Then $g(bt/a)$ has the leftmost point
of its support at $a$. If $g$ is invariant under Fourier
then we use $\sqrt{\frac ba}g(bt/a) + \sqrt{\frac
ab}g(at/b)$ to obtain again an invariant function, with
leftmost point of its support at $a$.
\end{proof}

\begin{definition}\label{def:Lproperty}
We say that a function $F(s)$, analytic
  in $\CC$ with at most  finitely many poles, has the
  \emph{L-Property} if the estimates $F(\sigma+i\tau) =
  O_{a,b,\epsilon}((1+|\tau|)^{(\frac12 - a)^++\epsilon})$
  hold (away from the poles), for $-\infty<a\leq\sigma\leq
  b<\infty$, $\epsilon>0$.
\end{definition}

\begin{theorem}\label{lindelof} The functions in
$\wh{L_a}$ have the L-Property. \end{theorem}

\begin{proof} Let $g(t)$ be a function in $L_a$ and let $G(s) =
\int_0^\infty g(t)t^{-s}\,dt$ be its right Mellin
transform. The function $g(t)$ is  a constant $\alpha(g)$ on
$(0,a)$. An expression for $G(s)$ as a meromorphic function
(in $\frac12<\Re(s)<1$, hence) in the right half-plane is:
\begin{eqnarray*} G(s) &=&
\frac{-\alpha(g)a^{1-s}}{s-1} + \int_a^\infty
g(t)t^{-s}\,dt\\ &=& \frac{-\alpha(g)a^{1-s}}{s-1}  +
\int_0^\infty \cF_+(\Un_{t>a}t^{-s})(u)\cF_+(g)(u)\,du\\
&=&\frac{-\alpha(g)a^{1-s}}{s-1} + \alpha(\cF_+(g))\int_0^a
\cF_+(\Un_{t>a}t^{-s})(u)du\\&+& \int_a^\infty
\cF_+(\Un_{t>a}t^{-s})(u)\cF_+(g)(u)\,du\\ \end{eqnarray*}
We established in \cite{cras2} a few results of an
elementary nature about the functions
$\cF_+(\Un_{t>a}t^{-s})(u)$ which are denoted there
$C_a(u,1-s)$ (in particular we showed that these functions
are entire functions of $s$).  For $\Re(s)<1$ one has
according to \cite[eq. 1.3.]{cras2}:
$$\cF_+(\Un_{t>a}t^{-s})(u) = \chi(s)u^{s-1}
-2\sum_{j=0}^\infty \frac{(-1)^j}{(2j)!}(2\pi
u)^{2j}\frac{a^{2j+1-s}}{2j+1-s}$$ hence for $0<\Re(s)<1$: $$\int_0^a
\cF_+(\Un_{t>a}t^{-s})(u)du =
\frac{\chi(s)a^s}{s}-2\sum_{j=0}^\infty
\frac{(-1)^j(2\pi)^{2j}}{(2j)!}\frac{a^{4j+2-s}}{(2j+1)(2j+1-s)}$$
This
is bounded on $\frac14\leq\Re(s)\leq\frac34$ (using
the well-known uniform estimate $|\chi(s)| \sim
|\Im(s)/2\pi|^{-\Re(s)+1/2}$ as $|\Im(s)|\to\infty$ in
vertical strips \cite[IV.12.3.]{tit}). We also have from
integration by parts and analytic continuation to $\Re(s)>0$
the expression:
$$\cF_+(\Un_{t>a}t^{-s})(u) = \frac{s\int_a^\infty
\sin(2\pi\,ut)t^{-s-1}\,dt - a^{-s}\sin(2\pi ua)}{\pi u}$$
which is $O(|s|/u)$ on $\frac14\leq\Re(s)\leq\frac34$,
$0<u$. Combining all this we find the estimate: $$G(s) =
O(|s|){\rm\ on\ } \frac14\leq\Re(s)\leq\frac34$$ This
(temporary) estimate justifies the use of the
Phragm\'en-Lindel\"of principle from bounds on $\Re(s) =
\frac12\pm\epsilon$. On any half-plane
$\Re(s)\geq\frac12+\epsilon>\frac12$ (excluding of course a
neighborhood of $s=1$) one has  $G(s) =
O(A^{\Re(s)})$ from the fact that $(s-1)G(s)/s$ belongs
to $A^s\HH^2$ and that elements of $\HH^2$ are bounded in
$\Re(s)\geq\frac12+\epsilon>\frac12$. And the functional
equation $G(1-s)=\chi(1-s)\wh{\cF_+(g)}(s)$ gives us
estimates on the left half-plane.  This shows that the
L-Property holds for $G(s)$. In particular, the Lindel\"of
exponents $\mu_G(\sigma)$ are at most $0$ for
$\sigma\geq\frac12$ and at most $\frac12 - \sigma$ for
$\sigma\leq\frac12$.
\end{proof}

\begin{remark}
In fact, the proof given above establishes the L-Property
for $G(s)$ in a stronger form than stated in the definition
\ref{def:Lproperty}. One has for example  $G(s) =
O_{\eta}\big((1+|\Im(s)|)^\eta)$ for each $\eta>0$, on the
strip $\frac12-\eta\leq \Re(s)\leq\frac12$, $G(s) =
O_\epsilon(|s|^\epsilon)$ for each $\epsilon>0$ on
$\frac12\leq \Re(s) \leq 1$ (away from the allowed pole at
$s=1$), and $G(s) = O_\eta(A^{\Re(s)})$ on
$\Re(s)\geq\frac12+\eta$, $\eta>0$.
\end{remark}

\begin{definition}
We let ${\mathcal L}_1$ to be the sub-vector space of $L_1$
containing the functions $g(t)$ whose right-Mellin transforms
$G(s)$ are $O_{g, a,b,N}(|s|^{-N})$ on all vertical strips $a\leq
\Re(s)\leq b$, and for all integers $N\geq1$  (away from the
pole, and the implied constant depending on $g$, $a$, $b$, and
$N$).
\end{definition}

\begin{theorem}\label{densethm} The sub-vector space ${\mathcal L}_1$
is dense in $L_1$. \end{theorem}

\begin{proof} \From\  proposition \ref{denselemma} we only
have to show that any function $G(s)$ in a $\wh{L_b}$, $b>1$
is in the closure of $\wh{{\mathcal L}_1}$. For this let
$\theta(s)$ be the Mellin transform of a smooth function
with support in $[1/e,e]$, satisfying $\theta(\frac12) =
1$. The function $\theta(s)$ is an entire function which
decreases faster than any (inverse) power of $|s|$ as
$|\Im(s)|\to\infty$ in any given strip $a\leq\sigma\leq
b$. Let us consider the functions $G_\epsilon(s) =
\theta(\epsilon (s-\frac12)+\frac12)G(s)$ as $\epsilon\to
0$. On the critical line they are dominated by a constant
multiple of $|G(s)|$ so they are square-integrable and
converge in $L^2$-norm to $G(s)$. We prove that for
$1\leq\exp(-\epsilon)b<b$ these functions all belong to
${\mathcal L}_1$. Their quick decrease in vertical strips is
guaranteed by the fact that $G(s)$ has the L-Property.  The
function
$$\frac{s-1}s G_\epsilon(s) = \theta\big(\epsilon
(s-\frac12)+\frac12\big)\frac{s-1}s G(s)$$ on the critical
line is the Mellin transform of a multiplicative convolution
on $(0,\infty)$ of an element in $L^2(b,\infty)$ with a
smooth function supported in
$[\exp(-\epsilon),\exp(+\epsilon)]$. The support of this
multiplicative convolution will be included in $[1,\infty)$
if $1\leq\exp(-\epsilon)b$. So for those $\epsilon>0$  one
has $G_\epsilon(s)\in \frac{s}{s-1}\HH^2$. Its image under
$\cF_+$ is $\theta(-\epsilon (s-\frac12)+\frac12)\cF_+(G)(s)
= \theta^\tau(\epsilon (s-\frac12)+\frac12)\cF_+(G)(s)$
where $\theta^\tau(w) = \theta(1-w)$ has the same properties
as $\theta(w)$ (we recall our abusive notation $\cF_+(F)(s)
= \chi(s)F(1-s)$.) Hence $\cF_+(G_\epsilon)(s)$ also belongs
to $\frac{s}{s-1}\HH^2$ and this completes the proof that
$G_\epsilon(s)\in \wh{{\mathcal L}_1}$.
\end{proof}

\begin{lemma} The subspace ${\mathcal L}_1$ is stable under $\cF_+$.
\end{lemma}

\begin{proof} Clear from the estimates of $\chi(s)$ in
vertical strips (\cite[IV.12.3.]{tit}). \end{proof}

\section{Completeness of the system of functions $\zeta(s)/(s-\rho)$}

We will use a classical estimate on the size of
$\zeta(s)^{-1}$:

\begin{proposition}[from {\cite[IX.7.]{tit}}]\label{inversezeta}
There is a real number $A$ and a strictly increasing
sequence $T_n>n$ such that $|\zeta(s)|^{-1} < |s|^A$ on
$|\Im(s)|=T_n$, $-1\leq \Re(s) \leq +2$. \end{proposition}

\begin{note}\label{note:sommezeros}
\From\  now on an infinite sum $\sum_\rho a(\rho)$ (with
complex numbers or functions or Hilbert space vectors
$a(\rho)$'s indexed by the non-trivial zeros of the Riemann
zeta function) means $\lim_{n\to\infty} \sum_{|\Im(\rho)|<
T_n} a(\rho)$, where the limit might be,  if we are dealing
with functions, a pointwise almost everywhere limit, or a
Hilbert space limit. When we say that the partial sums are
bounded (as complex numbers, or as Hilbert space vectors) we
only refer to the partial sums as written above. When we say
that the series is absolutely convergent it means that we
group together the contributions of the $\rho$'s with $T_n <
|\Im(\rho)| < T_{n+1}$ before evaluating the absolute value
or Hilbert norm. When building series of residues we write
sometimes things as if the zeros were all simple: this is
just to make the notation easier, but no hypothesis is made
in this paper on the multiplicities $m_\rho$, and the
formula used for writing $a(\rho)$ is a symbolic
representation, valid for a simple zero, of the more
complicated expression which would apply in case of
multiplicity, which we do not spell out explicitely.
\end{note}

\begin{theorem}\label{main} Let $G(s)$ be a function in
$\wh{L_1}$ which belongs to the dense subspace
$\wh{{\mathcal L}_1}$ of functions with quick decrease in
vertical strips. Then the series of residues for a fixed
$Z\neq1$, not a zero:
$$\sum_\rho
\frac{G(\rho)}{\zeta^\prime(\rho)}\frac{\zeta(Z)}{Z -
\rho}$$ converges absolutely pointwise to $G(Z)$ on
$\CC\setminus\{1\}$. It also converges absolutely in
$L^2$-norm to $G(Z)$ on the critical line. \end{theorem}

This is a series of residues for
$\frac{G(s)}{\zeta(s)}\frac{\zeta(Z)}{Z-s}$ where $s$ is the
variable and $Z\neq1$ is a parameter (with the exception of
the residue at $s=Z$). We have written the contribution of
$\rho$ as if it was simple (Titchmarsh uses a simlar
convention in \cite[IX.8.]{tit}). In fact the exact
expression is a linear combination of $\zeta(Z)/(Z -
\rho)^l$, $1\leq l \leq m_\rho$. We note that the trivial
zeros and $s=1$ are not singularities and contribute no
residue. 

\begin{proof}[Proof of Theorem \ref{main}] Let us consider
first the pointwise convergence. We fix $Z$, not $1$ and not
a zero and consider the function of $s$
$$\frac{G(s)}{\zeta(s)}\frac{\zeta(Z)}{Z-s}$$
We apply the
calculus of residues to the contour integral around a
rectangle with corners $\frac12\pm A\pm iT_n$ where $A\geq
\frac32$ is chosen sufficiently large such that both $Z$ and
$1-Z$ are in the open rectangle when $n$ is large enough.
Thanks to \ref{inversezeta} and the fact that $G(s)$ has
quick decrease the contribution of the horizontal segments
vanish as $n\to\infty$. The contribution of the vertical
segments converge to the (Lebesgue convergent) integral over
the vertical lines and we obtain:
$$\sum_{\rho}
\frac{G(\rho)}{\zeta^\prime(\rho)}\frac{\zeta(Z)}{Z - \rho}
- G(Z) = \frac1{2\pi}\left(\int_{\Re(s) = \frac12 + A} -
\int_{\Re(s) = \frac12 -
A}\right)\frac{\zeta(Z)G(s)}{(Z-s)\zeta(s)}|ds|$$
We prove
that the vertical contributions vanish. The functional
equation
$$\frac{G(1-s)}{\zeta(1-s)} =
\frac{\cF_+(G)(s)}{\zeta(s)}$$
reduces the case $\Re(s) =
\frac12 - A$ to the case $\Re(s) = \frac12 + A$. That last
integral does not change when we increase $A$. We note that
the $L^2$-norms of $G(s)$ on $\Re(s)=\sigma\geq2$ are
uniformly bounded because this is true with $G(s)$ replaced
with $(s-1)G(s)/s$ (which belongs to a Hardy space). Also
$1/\zeta(s) = O(1)$ in $\Re(s)\geq2$. The Cauchy-Schwarz
inequality then shows that the integral goes to $0$ as
$A\to\infty$. This proves the pointwise convergence:
$$G(Z) = \sum_{\rho}
\frac{G(\rho)}{\zeta^\prime(\rho)}\frac{\zeta(Z)}{Z -
\rho}$$ Going back to the contribution of the zeros with
$T_n<|\Im(\rho)|<T_{n+1}$, and expressing it as a contour
integral we see using $G(s)\in {\mathcal L}_1$ and
Proposition \ref{inversezeta} that the series of residues is
clearly absolutely convergent (with the meaning explained in
Note \ref{note:sommezeros}).

To show that the series converges to $G(Z)$  in $L^2$ on the
critical line it will be enough to prove it to be absolutely
convergent in $L^2$. We may with the same kind of reasoning
prove the absolute convergence of the series of residues:
$$\sum_{\rho} \frac{G(\rho)}{\zeta^\prime(\rho)}$$
For
this we consider $G(s)/\zeta(s)$ along rectangles with
vertical borders on $\Re(s) = \frac12\pm\frac32$ and
horizontal borders at the $\pm T_n$ and $\pm T_{n+1}$. The
functional equation (to go from $\Re(s)=-1$ to $\Re(s)=2$),
the estimate \ref{inversezeta} and the fact that $G(s)$
belongs to ${\mathcal L}_1$ then combine to prove that this series of
residues is absolutely convergent. In fact it converges to
$0$ as we prove later, but this is not needed here. So
returning to the problem of $L^2$-convergence we need only
prove the $L^2$ absolute convergence on the critical line
of:
\begin{eqnarray*} & &\sum_{\rho}
\frac{G(\rho)}{\zeta^\prime(\rho)}\left(\frac{\zeta(Z)}{Z -
\rho} - \frac{\zeta(Z)}{Z + 2}\right)\\ &=&\sum_{\rho}
\frac{G(\rho)}{\zeta^\prime(\rho)}\zeta(Z)\frac{2+\rho}{(Z -
\rho)(Z+2)} \\ \end{eqnarray*}
And for this it will be
sufficient to prove the $L^2$ absolute convergence of:
$$\sum_{\rho}
\frac{G(\rho)}{\zeta^\prime(\rho)}\zeta(Z)\frac{Z-1}{Z+2}\frac{2
+\rho}{(Z-\rho)(Z+2)}$$
We note that the function
$(Z-1)\zeta(Z)/(Z+2)^3$  belongs to $\HH^2(\Re(s)\geq
-\frac12)$, hence the same holds for each of the function
above depending on $\rho$. In case of a multiple zero its
contribution must be re-interpreted as a residue and will be
as a function of $Z$ a linear combination of the
$(Z-1)\zeta(Z)/(Z-\rho)^l(Z+2)^2$, $1\leq l \leq m_\rho$,
which also belong to $\HH^2(\Re(s)\geq -\frac12)$. It will
thus be enough to prove that the series above is
$L^2$-absolutely convergent on the line $\Re(Z) = -\frac12$,
\emph{as the norms are bigger on this line than on the
critical line}. We may then remove one factor $(Z-1)/(Z+2)$
and we are reduced to show that
$$\sum_{\rho}
\frac{G(\rho)}{\zeta^\prime(\rho)}\frac{\zeta(Z)(2
+\rho)}{(Z-\rho)(Z+2)}$$
is $L^2$-absolutely convergent on
$\Re(Z) = -\frac12$. What we do now is to reexpress for each
$Z$ on this line the contributions of the zeros with $T_n<
|\Im(\rho)| <T_{n+1}$ as a contour integral on the
rectangles (one with positive imaginary parts and the other
its reflection in the horizontal axis) bordered vertically
by $\Re(s) = -\frac14$ and $\Re(s) = +\frac54$. This will
involve along this contour the function of $s$:
$$\frac{G(s)}{\zeta(s)}\frac{\zeta(Z)(2+s)}{(Z-s)(Z+2)}$$
For a given fixed $s$ with $-\frac14\leq \Re(s)\leq+\frac54$
the function of $Z$ on $\Re(Z) = -\frac12$ given by
$$\frac{\zeta(Z)(2+s)}{(Z+2)(Z-s)}$$
has its $L^2$ norm
which is $O(1+|s|^2)$. Indeed:
$$\frac{2+s}{Z-s}
=\frac{2+s}Z\frac{Z}{Z-s}= \frac{2+s}Z (1 + \frac{s}{Z-s}) =
\frac{O(1+|s|^2)}{|Z|}$$
and $\zeta(Z)/Z(Z+2)$ is
square-integrable on $\Re(Z) = -\frac12$. The integrals
along these rectangular contours of the absolute values
$(1+|s|^2) |G(s)|/|\zeta(s)|$ give, from the quick decrease
of $G(s)$ and the Proposition \ref{inversezeta}, a convergent
series. With this the proof of \ref{main} is complete.
\end{proof}

This gives:

\begin{corollary}\label{cor:main}
 The functions $\zeta(s)/(s - \rho)^l$,
$1\leq l \leq m_\rho$ associated with the non-trivial zeros
are a complete system in $\wh{L_1}$. \end{corollary}

We also take note of the following:

\begin{proposition}\label{sumthm} 
One has for each $G(s)$ in the dense subspace $\wh{{\mathcal
L}_1}$:
$$0 = \sum_{\rho}
\frac{G(\rho)}{\zeta^\prime(\rho)}$$
where the series of
residues is absolutely convergent. 
\end{proposition}

\begin{proof} We have indicated in the proof of \ref{main}
that the series is absolutely convergent and its value is
$$\frac1{2\pi}\left(\int_{\sigma = 2} - \int_{\sigma =
-1}\right)\frac{G(s)}{\zeta(s)}|ds|$$
We prove that the
$\sigma = 2$ integral vanishes, and the $\sigma=-1$ integral
will then too also from the functional equation
$$\frac{G(1-s)}{\zeta(1-s)} = \frac{\cF_+(G)(s)}{\zeta(s)}$$
where $\cF_+(G)$ also belongs to $\wh{{\mathcal
L}_1}$. Using on $\sigma=2$ the absolutely convergent
expression $\frac1{\zeta(s)} = \sum_{k\geq1} \mu(k) k^{-s}$
it will be enough to prove:
$$0 = \frac1{2\pi}\int_{\sigma = 2}G(s)
k^{-s}|ds|$$
We shift the integral to the critical line
and obtain
$$\frac1{2\pi}\int_{\sigma = \frac12}G(s) k^{-s}|ds| +
\frac{{\rm Res}_1(G)}k$$ On the critical line we have in the
$L^2$-sense $G(s) = \int_0^\infty f(t)t^{s-1}dt$ for a
certain square-integrable function $f(t)$ (which is
$g(1/t)/t$ with $G(s) = \wh{g}(s)$). The Fourier-Mellin
inversion formula gives, in square-mean sense:
$$f(t) =
\lim_{T\to+\infty}\frac1{2\pi}\int_{s=\frac12-iT}^{s=\frac12+iT}G(s)
t^{-s}|ds|$$ As $G(s)$ is $O(|s|^{-N})$ on the critical line
for arbitrary $N$, we find that $f(t)$ is a smooth function
on $(0,\infty)$ given pointwise by the above formula. \From\  the
definition of $L_1$ one has $f(t) = c/t$ for $t\geq1$ with
a certain constant $c$.  The function $\int_0^1
f(t)t^{s-1}dt$ is analytic for $\Re(s)>\frac12$ so the
residue of $G(s)$ comes from $\int_1^\infty c\cdot
t^{s-2}dt$. This is first for $\Re(s)<1$ then by analytic
continuation the function $-c/(s-1)$ so ${\rm Res}_1(G) =
-c$, and on the other hand $f(k) = +c/k$. Combining all this
information the proof is complete. \end{proof}

\begin{remark} The result is (slightly) surprising at first
as we will prove that the evaluators associated with the
zeros are a complete and minimal system. \end{remark}

\begin{remark} These computations of residues are
reminiscent of a formula of Ramanujan which is mentioned in
Titchmarsh \cite[IX.8.]{tit}. For $ab=\pi, a>0$:
$$%
 \sqrt{a}\sum_{n=1}^\infty\frac{\mu(n)}n e^{-(a/n)^2} -
 \sqrt{b}\sum_{n=1}^\infty\frac{\mu(n)}n e^{-(b/n)^2} =
 -\frac1{2\sqrt b}\sum_\rho
 b^\rho\frac{\Gamma(\frac{1-\rho}2)}{\zeta^\prime(\rho)}
$$ where the meaning of the sum over the zeros is the one
from Note \ref{note:sommezeros}. \end{remark}

\section{Completion of the proofs of \ref{thmA}, \ref{thmAprime}, \ref{thmB}}

We also prove that the evaluators associated with the zeros
are a complete system in $L_1$. This is a further
application of the technique of \cite[Chap. 6]{hab} which
uses the theory of Nevanlinna functions and especially that
part of a fundamental theorem of Krein \cite{krein}
which says that an entire function which is Nevanlinna in
two complementary half-planes is necessarily of finite
exponential type (see {e.g.} \cite[I.\S4]{gorba}).

\begin{proposition} Let $a\geq1$.
 The vectors $Y^a_{\rho,k}$ associated
with the non-trivial zeros of the Riemann zeta function are
complete in $L_a$. 
\end{proposition}

\begin{proof} If $g\in L_a$ is
perpendicular to all those vectors (hence also $\cF_+(g)$)
then its right Mellin transform $G(s)$ factorizes as:
$$G(s) = \zeta(s)\theta(s)$$ with an entire function
$\theta(s)$. We have used that $G(s)$ shares with $\zeta(s)$
its trivial zeros and has at most a pole of order $1$ at
$s=1$. This expression proves that $\theta(s)$ belongs to
the Nevanlinna class of the right half-plane (as $G(s)$ and
$\zeta(s)$ are meromorphic functions in this class). \From\
the functional equation:
$$\theta(1-s) = \frac{G(1-s)}{\zeta(1-s)} =
\frac{\wh{\cF_+(g)}(s)}{\zeta(s)}$$
we see that $\theta(s)$ also belongs to the
Nevanlinna class of the left half-plane. According to the
theorem of Krein \cite{krein,gorba} it is of finite
exponential type which (if $\theta$ is not the zero
function) is given by the formula:
$$\max(\limsup_{\sigma\to+\infty}
\frac{\log|\theta(\sigma)|}{\sigma},
\limsup_{\sigma\to+\infty}
\frac{\log|\theta(1-\sigma)|}{\sigma})$$ We know that
$G(s)/\zeta(s)$ is $O(A^{\Re(s)})$ (with $A= 1/a$) in
$\Re(s)\geq2$ and similarly for $\cF_+(G)(s)/\zeta(s)$. So
this settles the matter for $A<1$ ($a>1$) as the formula
gives a strictly negative result. For $a=1$ we obtain that
$\theta(s)$ is of minimal exponential type. \From\  the
expression $G(s)/\zeta(s)$ on $\Re(s) = 2$ it is
square-integrable on this line. \From\  the Paley-Wiener
Theorem \cite{palwin} being of minimal exponential type it
in fact vanishes identically.
\end{proof}

\begin{proposition} 
Let $a>1$.  The vectors $Y^a_{\rho,k}$ associated with the
non-trivial zeros of the Riemann zeta function are not
minimal: indeed they remain a complete system in $L_1$ even
after omitting arbitrarily finitely many among them.
\end{proposition}

\begin{proof}
We adapt the proof of the preceding proposition to omitting the
vectors associated with the zeros from a finite set $R$. The
starting point will be
$$G(s) = \frac{\zeta(s)}{\prod_{\rho\in R}
(s-\rho)^{m_\rho}}\theta(s)$$ for a certain entire function
$\theta(s)$. The Krein formula for its exponential type
again gives a strictly negative result. So $\theta$ vanishes
identically.
\end{proof}

\begin{theorem}
Let $a=1$. The vectors $Y^1_{\rho,k}$ associated with the
non-trivial zeros of the Riemann zeta function are a minimal
(and complete) system
in $L_1$. The vectors, inverse Mellin transforms of the
functions $\zeta(s)/(s-\rho)^l$,  $1\leq l \leq m_\rho$, are a
minimal (and complete) system in $L_1$.
\end{theorem}

\begin{proof}
The fact that the functions $\zeta(s)/(s-\rho)^l$, $1\leq l
\leq m_\rho$ belong to $\wh{L_1}$ implies that the
evaluators $Y^1_{\rho,k}$'s are a minimal system. We know
already that they are a complete system. The system of the
$\zeta(s)/(s-\rho)^l$, $1\leq l \leq m_\rho$, is, up to
triangular invertible linear combinations for each $\rho$
the uniquely determined dual system. As a dual system it has
to be minimal. And we know already from \ref{cor:main} that
it is a complete system.
\end{proof}

This completes the proof of \ref{thmB}.

\begin{proposition}
Let $a<1$. The vectors $Y^a_{\rho,k}$ are minimal and not
complete in $L_a$.
\end{proposition}

\begin{proof}
If they were not minimal, their orthogonal projections to
$L_1$ which are the vectors $Y^1_{\rho,k}$, would not be
either. And they are not complete from the existence of the
co-Poisson subspace $P_a$.
\end{proof}

With this the proof of \ref{thmA} is completed, with the
exception of the identification of the co-Poisson space as
the perpendicular complement to the space spanned by the
$Y^a_{\rho,k}$'s. We refer the reader to \cite[Chap.6]{hab}
especially to \cite[Theorems 6.24, 6.25]{hab} which have all
the elements for the proof, as it does not appear useful to
devote space to this here.

\begin{proposition} The vectors $Z^1_{\rho,k}$ are not minimal
in $K_1$. In fact $K_1$ is spanned by these vectors even
after omitting $Z^1_{\rho_1,m_{\rho_1}-1}$ and $Z^1_
{\rho_2,m_{{\rho_2}}-1}$ ($\rho_1\neq\rho_2$), or
$Z^1_{\rho,m_\rho-1}$ and $Z^1_ {\rho,m_{\rho}-2}$
($m_\rho\geq2$), from the list. This shortened system is
then a minimal system.
\end{proposition}

\begin{proof} If $f$ in $K_1$ is perpendicular (for the
  form $[f,g] = \int_1^\infty f(t)g(t)dt$) to this
  shortened list of evaluators then its right Mellin
  transform  factorizes as
$$F(s) =
\frac{s(s-1)\zeta(s)}{(s-\rho_1)(s-\rho_2)}\theta(s)$$
where we have used that $F(0) = 0$ and that $F(s)$ has no
pole at $s=1$. In this expression we have the two cases
$\rho_1\neq\rho_2$ and $\rho_1 = \rho_2$. The proof
then proceeds as above and leads to $F(s) = 0$. 
To prove minimality for the shortened system one only has to
consider the functions
$$\frac{s(s-1)\zeta(s)}{(s-\rho_1)(s-\rho_2)}\frac1{(s-\rho)^l}$$
associated with the remaining zeros (and remaining
multiplicities), as they are easily seen to be the right Mellin
transforms of elements from the Sonine space $K_1$.
\end{proof}

\begin{proposition}
Let $a>1$. The vectors  $Z^a_{\rho,k}$ span $K_a$ even after
omitting arbitrarily finitely many among them.
\end{proposition}

\begin{proof}
They are the orthogonal projections to $K_a$ of the vectors
$Y^a_{\rho,k}$ in $L_a$.
\end{proof}

\begin{theorem}
Let $a<1$.  The vectors $Z^a_{\rho,k}$ are minimal
in $K_a$.
\end{theorem}

\begin{proof}
Let $\theta(t)$ be a smooth non-zero function supported in
$[a,A]$ ($A= 1/a>1$). Its right Mellin transform
$\wh{\theta}(s)$ is then $O(A^{|\Re(s)|})$ on $\CC$. And if
$P(s)$ is an arbitrary polynomial, then $P(s)\wh{\theta}(s)
= \wh{\theta_P}(s)$ for a certain smooth function
$\theta_P$, again supported in $[a,A]$, so
$\wh{\theta_P}(s)= O_P(A^{|\Re(s)|})$. Hence
$\wh{\theta}(s)$ decreases faster than any inverse
polynomial in any given vertical strip, in particular on
$-1\leq\Re(s)\leq2$. \From\ this we see that the function
$G(s)= s(s-1)\wh{\theta}(s)\zeta(s)$ is square-integrable on
the critical line and belongs to $A^s\HH^2$ (one may write
$G(s) = s^3\wh{\theta}(s)(s-1)\zeta(s)/s^2$, and use the
fact that $(s-1)\zeta(s)/s^2$ belongs to $\HH^2$). We have
$\chi(s)G(1-s) = s(s-1)\wh{\theta}(1-s)\zeta(s)$ so again
this belongs to $A^s\HH^2$. This means that $G(s)$ is the
right Mellin transform of a (non-zero) element $g$ of
$K_a$. Let us now take a non-trivial zero $\rho$, which for
simplicity we assume simple. We choose the function
$\theta(t)$ to be such that $\wh{\theta}(\rho)\neq0$, which
obviously may always be arranged. Then, using
\ref{zetaprop},  $G(s)/(s-\rho)$ is again the Mellin
transform of a non-zero element $g_\rho$ in $K_a$. This
element is perpendicular (for the bilinear form $[f,g]$) to
all the evaluators except $Z^a_{\rho,0}$, to which it is not
perpendicular. So $Z^a_{\rho,0}$ can not be in the closed
span of the others. The proof is easily extended to the case
of a multiple zero (we don't do this here, as the next
section contains a proof of a more general statement).
\end{proof}

The three theorems \ref{thmA}, \ref{thmAprime}, \ref{thmB}
are thus established.

\section{Zeros and evaluators for general Sonine functions}\label{sec:zeros}

Let us more generally associate to any non-empty multiset
$\cZ$ of complex numbers (a countable collection of complex
numbers, each assigned a finite multiplicity) the problem of
determining whether the associated evaluators are minimal,
or complete in a Sonine space $K_a$ or an extended Sonine
space $L_a$. To be specific we consider the situation in
$K_a$, the discussion could be easily adapted to
$L_a$. \From\ 
the fact that the Sonine spaces are a decreasing chain, with
evaluators in $K_a$ projecting orthogonally to the
evaluators in $K_b$ for $b\geq a$, we may associate in
$[0,+\infty]$ two indices $a_1(\cZ)$ and $a_2(\cZ)$ to the
multiset $\cZ\in\CC$. The index $a_1(\cZ)$ will be such that
the evaluators are a minimal system for $a<a_1(\cZ)$ and not
a minimal system for $a>a_1(\cZ)$ and the index $a_2(\cZ)$
will be such that the evaluators are complete for
$a>a_2(\cZ)$ but not complete for $a<a_2(\cZ)$. Let us take
for example the multiset $\cZ$ to have an accumulation point
$w$ (there is for each $\epsilon>0$ at least one complex
number $z$ in the support of $\cZ$ with $0<|z-w|<\epsilon$):
then the system is never minimal and is always complete so
that $a_1 = 0$ and $a_2 = 0$. As another example we take the
multiset to have finite cardinality: then the evaluators are
always minimal and never complete so $a_1 = +\infty$, and
$a_2 = +\infty$. For the zeros of the Riemann zeta function
we have $a_1 = a_2 = 1$. There is a general phenomenon here:

\begin{theorem}\label{thmC}
The equality $a_1(\cZ) = a_2(\cZ)$  always holds.
\end{theorem}

Let us thus write $a(\cZ)$ for either $a_1(\cZ)$ or
$a_2(\cZ)$. We will prove that $a(\cZ)$ does not change from
adding or removing a finite multiset to $\cZ$ (maintaining
$\cZ$ non-empty):

\begin{theorem}\label{thmD1}
If $0<a<a(\cZ)$ then the evaluators associated to $\cZ$ remain
not complete, and minimal, in $K_a$, after including
arbitrarily finitely many other evaluators.
\end{theorem}

\begin{theorem}\label{thmD}
If $a(\cZ)<a<\infty$ then the evaluators associated to $\cZ$
remain  complete, and not minimal,  in $K_a$ after omitting
arbitrarily finitely many among them.
\end{theorem}

We will say that $g(t)$ is a Sonine function if it belongs
to $\cup_{a>0} K_a \subset L^2(0,\infty;dt)$. We also say
that $G(s)$ is a Sonine function if it is the right Mellin
transform of such a $g(t)$.

\begin{lemma}\label{lem:a2lessa1}
If the system of evaluators associated in a given $K_a$ to a
(non-empty) multiset $\cZ$ is not complete, then it is
minimal. Alternatively, if it is not minimal, it has to be complete.
\end{lemma}

\begin{proof}
Let us assume that the system is not complete. Then we have
a non zero Sonine function $G(s)$ in $\wh{K_a}$ such that
$\pi^{-s/2}\Gamma(\frac s2)G(s)$ vanishes on $\cZ$. \From\ the
proposition \ref{zeroprop} we know that if
$\pi^{-w/2}\Gamma(\frac w2)G(w) = 0$ then $G(s)/(s-w)$ is
again a Sonine function in $\wh{K_a}$. Let us now proceed to
take $\rho$ in the support of $\cZ$, and divide $G(s)$ by
powers of $(s-\rho)$ to construct functions which vanish
exactly to the $k$-th order at $\rho$, for $0\leq k <
m_\cZ(\rho)$ (this is after incorporating the Gamma
factor). \From\ suitable linear combinations we construct
further an $a$-Sonine function $G_k(s)$ whose $l$-th
derivative  for $0\leq l <m_\cZ(\rho)$ (again with the Gamma
factor incorporated) vanishes at $\rho$, except for $l=k$
for which it does not vanish, and with $G_k(s)$ vanishing on
the remaining part of the multiset $\cZ$. This proves that
the evaluators in $K_a$ associated with $\cZ$ are minimal.
\end{proof}

We note that this provides an alternative route to our
statement from \cite{cras2} that finitely many evaluators
are always linearly independent in $K_a$, once we know that
$K_a$ is infinite dimensional.

\begin{lemma}\label{lem:a1lessa2}
If the system of evaluators associated in a given $K_a$ to a
(non-empty) multiset $\cZ$ is minimal, then it is not
complete in any $K_b$ with $b<a$.
\end{lemma}

\begin{proof}
We pick a $\rho$ in the support of $\cZ$, with multiplicity
$m_\rho$. As the system is minimal, we have the existence of
at least one Sonine function $G(s)$ in $\wh{K_a}$ which
vanishes on the other part of $\cZ$ but vanishes only to the
$(m_\rho-1)$-th order at $\rho$. Let us now consider a
function $F(s) = \theta(s)G(s)$ where $\theta(s)$ is the
Mellin transform of a non-zero smooth function supported in
an interval $[\exp(-\epsilon), \exp(+\epsilon)]$. We know
from Theorem \ref{lindelof} that Sonine functions have the
L-Property, so using the arguments of the smoothing
technique in the proof of Theorem \ref{densethm} we obtain
easily that any such $F(s)$ is a non-zero element of
$\wh{L_{b}}$ for any $b \leq\exp(-\epsilon)a$. Replacing
$\theta(s)$ by $(s-\rho)\theta(s)$ we may impose
$\theta(\rho) = 0$. Then $F(s)$ (with the Gamma factor)
vanishes on $\cZ$ and this proves that the evaluators
associated with $\cZ$ are not complete in $\wh{L_{b}}$.
\end{proof}

At this stage we have completed the proof of Theorem
\ref{thmC}: Lemma \ref{lem:a2lessa1} implies $a_2(\cZ)\leq
a_1(\cZ)$ and Lemma \ref{lem:a1lessa2} implies $a_1(\cZ)
\leq a_2(\cZ)$.

\begin{lemma}
If the system of evaluators associated in a given $K_a$ to a
(non-empty) multiset $\cZ$ is minimal, then it is not
complete in any $K_b$ with $b<a$, even after adding to
the system of evaluators associated with $\cZ$ 
arbitrarily finitely many other evaluators.
\end{lemma}

\begin{proof}
We only have to replace the function $\theta(s)$ from the
preceding proof by $P(s)\theta(s)$ where $P(s)$ is an
arbitrary polynomial.
\end{proof}

This, together with Lemma \ref{lem:a2lessa1}, clearly implies
Theorem \ref{thmD1}. It also implies the Theorem \ref{thmD}:
let us suppose $a(\cZ)<a<\infty$. Let us imagine that after
removing finitely many evaluators we do not have a complete
system. Then this remaining system, being not complete, has
to be minimal from  Lemma \ref{lem:a2lessa1}. We just proved
that in these circumstances the system in a $K_b$ with $b<a$
can not be complete, even after including finitely many
arbitrary evaluators. This gives a contradiction for $a(\cZ)
< b < a$, as we may reintegrate the omitted evaluators. So
Theorem \ref{thmD} holds.

Let $g$ be a non-zero Sonine
function. We write $\lambda(g)>0$ for
the minimal point of the support of $g$ and $\mu(g)>0$ for
the minimal point of the support of $\cF_+(g)$. And we let
$a(g)$ be $\sqrt{\lambda(g)\mu(g)}$.

To each non-zero Sonine function $g$ we associate the
multiset $\cZ_g$ (which will be proven to have infinite
cardinality) of its non-trivial zeros: these are the zeros
of the completed Mellin transform $\pi^{-s/2}\Gamma(\frac
s2)\wh g(s)$, so $0$, $-2$, \dots, might be among them but
they are counted with multiplicity one less than in $\wh
g(s)$.

Before proceeding further we need to recall some classical
results from the Theory of Nevanlinna functions and Hardy
Spaces. We refer the reader for example to \cite[Chap.1]{bra} and
\cite[I.\S4]{gorba} for proofs and more detailed statements (see
also \cite{dym, dym2, hof}). A Nevanlinna function $F(s)$ in
a half-plane (we consider here $\Re(s)>\frac12$) is an
analytic function which may be written as the quotient of
two bounded analytic functions. To each non-zero $F$ is
associated a real-number $h(F)$, its mean type (in the
terminology from \cite{bra}), which may be obtained (in the
case of the half-plane $\Re(s)>\frac12$) from the formula
$h(F) = \limsup_{\sigma\to+\infty}
\log|F(\sigma)|/\sigma$. The mean-type of a product is the
sum of the mean types. The mean-type contributes a factor
$e^{h(s-\frac12)}$ to the Nevanlinna-Smirnov factorization
of the function $F(s)$, in particular, for the specific case
of the Smirnov-Beurling factorization of an element in
$\HH^2$, it gives the special inner factor (here
$h\leq0$). The other factors have mean type $0$. In the
particular case when $F(s)$ is the right Mellin transform of
a square-integrable function $f(t)$ supported in
$[\lambda,+\infty)$, $\lambda>0$, then the mean-type of $F$
is also $\log(\lambda(f)^{-1})$ where
$\lambda(f)\geq\lambda>0$ is the lowest point of the support
of $f$. To see this, we may after a multiplicative
translation assume that $\lambda(f) = 1$. We want to prove
that the mean-type of $F$ is $0$. One has $h\leq0$ as $F$ is
bounded, say for $\Re(s)\geq1$. In the canonical
factorization of $F$, the outer factor is still an element
of the Hardy space. The inner factor is bounded by $1$. So
$F(s)$ belongs to $e^{hs}\HH^2$, which is the subspace of
Mellin transforms of $L^2(e^{-h},\infty; dt)$, so $h=0$. We
conclude this brief summary with Krein's theorem
\cite{krein}, which we have already used in the
previous proofs. This important theorem (see \cite[I.\S4]{gorba})
states in particular
that an entire function $\theta(s)$ which is in the
Nevanlinna class in two complementary half-planes is
necessarily of finite exponential type. Furthermore the
exponential type is the maximum of the mean-types for the
two half-planes. Hence, if $\theta$ is not the zero
function, at least one of the two mean-types has to be
non-negative.

We will also need some classical results from the theory of
entire functions \cite{levin}. Let $F(z)$ be an entire
function. Then $F$ is said (\cite[I.\S12]{levin}) to have
normal type with respect to the (Lindel\"of) refined
(proximate) order $r\log(r)$ (which is the one useful to us
here) if
$$0 < \limsup_{r\to\infty} \frac{\log \max_{|z|=r}
|F(z)|}{r\log(r)} < \infty$$ If this holds, the generalized
Phragm\'en-Lindel\"of indicator function is defined as:
$$h_F(\theta) = \limsup_{r\to\infty}
\frac{\log|F(r\,e^{i\theta})|}{r\log(r)}$$ One proves that
the indicator function of the entire function $F(z)$ of
normal type is finite valued and is a continuous
``trigonometrically convex'' function of $\theta$
(\cite[I.\S18]{levin}).

\begin{remark}
The indicator function for $F(z-a)$ is the same as the one
for $F(z)$: to see this one may use the upper estimate
\cite[I.\S18, Thm 28]{levin} $\log|F(r\,e^{i\,\vartheta})|
< (h_F(\vartheta)+\epsilon)r\log(r)$ for $r>r_\epsilon$ in a
given open angular sector $|\arg(z) - \theta|<\eta$, and the
continuity of $h_F$ at $\theta$. The parallel ray starting
at $a$ is contained in this sector except for a finite
segment, so the indicator function based at $a$ is bounded
above by the one based at the origin, and vice versa.
\end{remark}

A ray $L_\theta = \{r\,e^{i\theta}, 0<r<\infty\}$ is a ray
of completely regular growth (CRG-ray, \cite[III]{levin})
for $F$ if
$$h_F(\theta) = \lim_{r\to\infty, r\notin E_\theta}
\frac{\log|F(r\,e^{i\theta})|}{r\log(r)}$$ where the
excluded set $E_\theta \subset (0,\infty)$ has vanishing
upper relative linear density. The set of CRG-rays is
closed. The entire function $F(z)$ is said to be of
completely regular growth if all the rays are CRG-rays. A
fundamental theorem \cite[III.\S3]{levin} which applies to
CRG-functions states that the number $n(r,\alpha,\beta)$ of
zeros of modulus at most $r$ in the open angular sector
$\alpha<\theta<\beta$ has the  following asymptotic behavior:
$$\lim_{r\to\infty} \frac{n(r,\alpha,\beta)}{r\log(r)} =
\frac1{2\pi}\Big(h^\prime_F(\beta) - h^\prime_F(\alpha) +
\int_\alpha^\beta h_F(\theta)d\theta\Big)$$ under the
condition that $h_F$ admits derivatives at  $\alpha$ and
$\beta$ (from the trigonometrical convexity 
right and left derivatives always exist).

\begin{theorem}\label{thm:zeros}
Let $g(t)$ be a non-zero Sonine function, with Mellin
transform $G(s)$, and Gamma-completed Mellin transform
$\cG(s)$. The entire function $\cG(s)$ is of normal type for
the Lindel\"of refined order $r\log(r)$. Its indicator
function is   $\frac12|\cos(\theta)|$. The entire function
$\cG(s)$ is a function of completely regular growth. The
number of its zeros of modulus at most $T$ in the angular
sector $|\arg(z-\frac12) - \frac\pi2|<\epsilon<\pi$ is
asymptotically equivalent to $\frac{T}{2\pi}\log(T)$, and
similarly for the angular sectors containing the lower-half
of the critical line. The number of zeros of $\cG(s)$ with
modulus at most $T$ in $|\arg(\pm z)|<\frac\pi2 -\epsilon$ is
$o(T)$.
\end{theorem}

\begin{proof}
We know from Theorem \ref{lindelof} that $G(s)$ has the
L-property, and in particular it is $O(1+|s|)$ in
$0\leq\Re(s)\leq1$, and furthermore it is $O(A^{\Re(s)})$ in
$\Re(s)\geq1$. On the other hand the Stirling formula easily
leads to
$$\lim_{r\to\infty} \frac{\max_{|s-\frac12|=r,\Re(s)\geq0}
\log |\pi^{-s/2}\Gamma(\frac s2)|}{r\log(r)} = \frac12$$ so
certainly
$$\limsup_{r\to\infty}
\frac{\max_{|s-\frac12|=r,\Re(s)\geq0} \log
|\cG(s)|}{r\log(r)} \leq \frac12$$ As further
$\limsup_{\sigma\to\infty} \log|G(\sigma)|/\sigma$ is finite
(it is the mean-type of $G$ in the right-half plane) we have
$\limsup_{\sigma\to\infty}
\log|G(\sigma)|/\sigma\log(\sigma) = 0$. Hence:
$$\limsup_{r\to\infty}
\frac{\max_{|s-\frac12|=r,\Re(s)\geq0} \log
|\cG(s)|}{r\log(r)} = \frac12$$ In $\Re(s)\leq1$ we have the
identical result as $\cG(1-s)$ is the completed Mellin
transform of $\cF_+(g)$.  So the entire function $\cG(s)$ is
of normal type $\frac12$ for the refined order
$r\log(r)$. The argument using the Stirling formula which
has led to the inequality above gives on any given ray with
$|\arg(\theta)|\leq\frac\pi2$ that its indicator function is
bounded above by $\frac12 \cos(\theta)$ (we use again that
for $|\theta|= \frac\pi2$ we have $G(s) = O(|s|)$, and for
$|\theta| < \frac\pi2$ we have $G(s) =
O_\theta(A^{\Re(s)})$). We show the property of complete
regular growth on a ray with $|\theta|<\frac\pi2$ and at the
same time identify the value of the indicator for this ray
to be $\frac12 \cos(\theta)$ (by continuity we will then
have the value of the indicator for
$|\theta|=\frac\pi2$). The Gamma factor gives the correct
limit and there are no excluded values for $r$; so we only
need to show that $\log|G(s)|/r\log(r)$ goes to $0$ as
$r\to\infty$ while avoiding an exceptional set
$E_\theta\subset (0,\infty)$ having vanishing upper relative
linear density ($r = |s-\frac12|$). Actually this holds with
$r$ replacing $r\log(r)$ and a finite (not necessarily zero)
limit, as $A^{-s}G(s)$ belongs to the Hardy space of the
right half-plane (and then we can 
invoke \cite[V.\S4, Thm 6]{levin}; 
it is all a matter of understanding the CRG-behavior
of a Blaschke product, as the other factors in the
canonical factorization are easily taken care
of). Going back to $\cG(s)$ we thus have its CRG property
(for the refined order $r\log(r)$) on the rays
$\arg(s-\frac12)=\theta$, $|\theta|<\pi/2$. Hence also in
the left half-plane as $\cG(1-s)$ is the completed Mellin
transform of $\cF_+(g)$ and the indicator function is thus
$\frac12|\cos(\theta)|$. As the set of CRG-rays is closed,
we conclude that $\cG(s)$ is a CRG entire function. The
central result from \cite[III]{levin} leads then to the
stated asymptotic densities of zeros in open sectors
containing either the upper half or the lower half of the
critical line. Concerning the sectors  $|\arg(\pm
z)|<\frac\pi2 -\epsilon$, the vanishing asymptotic linear
density of the zeros follows again from \cite[V.\S4,
Thm 6]{levin}, or more simply from the fact that $\sum
\frac1{|\rho|}$ converges for the zeros in such a sector
(the zeros of the Blaschke product satisfy $\sum
\frac{\Re(\rho) - \frac12}{|\rho|^2}<\infty$).
\end{proof}

\begin{remark}
In particular the function $\cG(s)$ is an entire function of
order one (and maximal type for this order) which admits a
representation as an  Hadamard product 
$s^Ne^{\alpha+\beta s}\prod_\rho (1-s/\rho)e^{s/\rho}$.
\end{remark}

\begin{remark}
The entire function $\cG(\frac12 + i z)$ is a function of
the ``class A'' as studied in \cite[V]{levin}.
\end{remark}

\begin{remark}\label{rem:tlogt}
In \cite{cras3} we have produced explicit formulae for some
even distributions $A_a(t)$ and $B_a(t)$ having the Sonine
property for the cosine transform.  We  proved that their
complete Mellin transforms $\cA_a(s)$ and $\cB_a(s)$ are the
structure functions of the de~Branges Sonine-cosine spaces
(no explicit formula had been known prior to \cite{cras3}):
and this has the interesting corollary that the Riemann
Hypothesis holds true for them.  We have presented in
\cite{cras4} a summary of further results of ours. The Note
contains formulae for some second order differential
operators intrinsically associated with the Fourier
Transform. Under suitable boundary conditions these
operators are self-adjoint with discrete spectrum, and the
squared imaginary parts of the zeros of $\cA_a(s)$ and
$\cB_a(s)$ are their eigenvalues (this proves in another
manner that $\cA_a(s)$ and $\cB_a(s)$ satisfy the Riemann
Hypothesis).  If $\rho$ is a fixed chosen zero then
$\cA_a(s)/(s-\rho)$ is the complete Mellin transform of a
square-integrable even function with the Sonine property. We
may then apply Theorem \ref{thm:zeros} with the result that
$\cA_a(s)$ (or $\cB_a(s)$)  share with
$\pi^{-s/2}\Gamma(\frac s2)\zeta(s)$ the principal order of
its asymptotic density of zeros.
\end{remark}

\begin{theorem}\label{thmE}
For each non-zero Sonine function there holds: $a(\cZ_g) =
a(g)$. This means in particular that $g$ has infinitely many
(non-trivial) zeros, that the evaluators associated to
$\cZ_g$ are minimal but not complete in $K_a$ if $a\leq
a(g)$ and that they are complete, even after omitting
arbitrarily finitely many among them, in $K_a$ if $a>a(g)$.
\end{theorem}

\begin{proof}
Let $g(t)$ be a non-zero Sonine function in $K_a$, with
Mellin transform $G(s)$. Let us consider the multi-set
$\cZ_g$ of the non-trivial zeros of $G(s)$ and the
associated evaluators in Sonine spaces $K_b$. Replacing $g$
by a multiplicative translate we may arrange that the lowest
point of its support coincides with the lowest point of the
support of $\cF_+(g)$, hence with the number we have denoted
$a(g)$, so we may assume $a=a(g)$. The system of evaluators
in $K_a$ is not complete, as all are perpendicular to $g$,
hence $a(\cZ_g)\geq a$. Let $b>a$ and let us prove that the
evaluators are complete in $K_b$. If not, there is a
non-zero function $f(t)$ in $K_b$ such that its Mellin
transform $F(s)$ factorizes as $F(s) = G(s)\theta(s)$ with
an entire function $\theta(s)$. In particular $\theta(s)$ is
a Nevanlinna function in the right half-plane. We know that
the mean types are related through $h(F) = h(G) +
h(\theta)$. We know that $h(G) = -\log a$ and that $h(F)
\leq -\log b$, hence $h(\theta) \leq \log(a/b) < 0$. On the
other hand we have $\chi(s)F(1-s) =
\chi(s)G(1-s)\theta(1-s)$. Repeating the argument for
$\theta(1-s)$ we obtain that its mean type is also
$<0$. According to Krein's theorem $\theta(s)$ has finite
exponential type given by the formula
$$\max(\limsup_{\sigma\to\infty}
\frac{\log|\theta(\sigma)|}{\sigma},
\limsup_{\sigma\to\infty}\frac{\log|\theta(1-\sigma)|}{\sigma}),
$$ hence we obtain a strictly negative result. This is
impossible, and the function $f\in K_b$ does not exist. So
the evaluators associated to $\cZ$ are complete in $K_b$ for
$b>a(g)$ and  $a(\cZ)\leq a = a(g)$. We know already
$a(\cZ)\geq a(g)$ so we have an equality, as was to be
proven. The other statements are just repetitions of
previously proven assertions.
\end{proof}

We extract from the proof above the following:

\begin{proposition}\label{thmF}
If  $f$ and $g$ are two non-zero Sonine functions such that
$\cZ_g \subset  \cZ_f$ then the entire function $F(s)/G(s)$
has finite exponential type and $a(f)\leq a(g)$.
\end{proposition}

\begin{proof}
We replace $g$ by a multiplicative translate so that $g\in
K_a$ with $a = a(g)$. And we similarly assume $f\in K_b$
with $b = a(f)$.  Krein's theorem is applied to the entire
function $\theta(s)$ with $F(s) = G(s)\theta(s)$, with the
conclusion that $\theta(s)$ has finite exponential
type. Furthermore exactly as in the previous proof if we had
$b>a$ we could prove that the mean types of $\theta(s)$ in
the left and right half-plane are both strictly negative,
which is impossible. So $b\leq a$, that is $a(f)\leq a(g)$.
\end{proof}

We also mention:

\begin{proposition}\label{thmG}
If  $f$ and $g$ are two non-zero Sonine functions such that
$\cZ_f = \cZ_g$ then $f$ and $g$ are multiplicative
translates of one another, up to multiplication by a
non-zero complex number.
\end{proposition}

\begin{proof}
The Hadamard product representation leads to an equality
$\wh f(s) = e^{\lambda + \mu s} \wh g(s)$. But Wiener's
theorem \cite{palwin} on the gain of a causal filter tells us that
$\log|\wh g(s)|$ and $\log|\wh f(s)|$ are both integrable against a
Cauchy weight on the critical line, and this implies that
$\mu$ has to be real. The equation then says exactly
that $f$ is, up to a multiplicative constant, a
multiplicative translate of $g$.
\end{proof}

\section{Conclusion}

Our theorems from
\cite[Chap. 6]{hab} concerning the completeness of the
evaluators associated to the Riemann zeta function and
the Dirichlet L-functions have been shown here
to be special instances of a more general
statement. Does this mean that these
theorems from \cite{hab} are not specific enough
to tell us anything
interesting? 

To discuss this, we shall, briefly,
mention a few basic aspects of the general theory of Hilbert
spaces of entire functions, and thus see why it is reasonable
to be hopeful of some connections with the problem of 
the Riemann hypothesis, at the
technical level at least (it is in the exact same manner,
no more no less, that the, more widely known, basic aspects of the Hilbert
theory of self-adjoint operators may be thought of bearing
some relevance to the technical aspects of the Riemann hypothesis).
We include this short paragraph despite the prolonged
existence of ethically unfortunate claims.
To each de~Branges
space are associated (up to some normalizations) a function
$\cA(z)$ and also a function $\cB(z)$ which both have all
their zeros on the symmetry axis. This is a corollary to the
way the functions $\cA(z)$ and $\cB(z)$ are related to the
Hilbert space structure, hence participates of the general
idea of thinking about the Riemann Hypothesis in Hilbert
space and operator-theoretical terms. We mention that the
work of de~Branges is closely related to the vast
investigations of M.G.~Krein \cite{gorba} on problems of
extrapolation of stationary processes, problems of
scattering theory, problems of moments, canonical systems,
\dots, where the operator theoretical aspects are quite
explicitely in the foreground. The zeta function is not
entire, but has only one pole. Its functional equation
involves Gamma factors, to which de~Branges associates the
two-parameter family of the Sonine spaces for the Hankel
transforms of parameter $\nu$, and the support conditions of
parameter $a$. We focus on the spaces associated with the
cosine and sine transforms. De~Branges \cite{bra92,bra94}
uses in his constructions the other Sonine spaces, even
``double-Sonine'' spaces: the idea of using $2$-dimensional
constructs to study a Riemann Hypothesis in dimension $1$ is
a familiar one from other contexts. The structure functions
$\cA_a(z)$ and $\cB_a(z)$, especially for the Sonine spaces
associated with the cosine transform, have analytic
properties and symmetries quite close to what is known to
hold for $\pi^{-s/2}\Gamma(s/2)\zeta(s)$, with one
interesting bonus: they are proven to satisfy the Riemann
Hypothesis. Explicit representations for these functions, as
completed Mellin transforms, have been obtained recently
(\cite{cras3}). We proved here (\ref{thm:zeros}) that they
have to first order the same density of zeros as the Riemann
zeta function. We have obtained (\cite{cras4}) a spectral
interpretation of their zeros, in terms of some Dirac and
Schr\"odinger operators which we have associated to the
Fourier Transform.

We explained in \cite{hab} the path which has led to our own
interest in all this: the path from the explicit formula
to the
co-Poisson formula and beyond. The co-Poisson formula leads to the
association with the zeta function of certain quotient
spaces of the Sonine spaces. We saw in the previous section
that some of the theorems originally proven for the zeta
function or the Dirichlet $L$-functions have more general
validity, as some aspects hold true for all Sonine
functions. The Riemann 
Hypothesis of course does not hold for all Sonine
functions (we may always add arbitrarily chosen zeros; 
it is also easy to construct an example of a
Sonine function with no zeros in the critical
strip). Nevertheless it might be that some other aspects,
known or expected to hold for the Riemann zeros, do have
some amount of wider validity; further investigations of the
zeros of Sonine functions are needed to better understand
the situation.


\begin{small}
\textbf{Acknowledgments.} I thank Michel Balazard
and \'Eric Saias, for discussion on Sonine spaces,
and especially on the functions $\zeta(s)/(s-\rho)$.\par
\end{small}


\bibliographystyle{amsplain}

\end{document}